\documentclass{birkjour}

 \theoremstyle{definition}
 
 \theoremstyle{remark}

 \numberwithin{equation}{section}

\begin{document}

\title[]
{On Resolving Singularities of Plane Curves via a Theorem  attributed to Clebsch}

\author[]{David E. Rowe\\Mainz University}


\email{rowe@mathematik.uni-mainz.de}


\maketitle

\noindent
Abstract: This paper discusses a central theorem in birational geometry first proved by Eugenio Bertini in 1891. J.L. Coolidge described the main ideas behind Bertini's proof, 
but he attributed the theorem to Clebsch. He did so owing to a short note that Felix Klein appended to the republication of Bertini's article in 1894. The precise circumstances that led to Klein's intervention can be easily reconstructed from letters Klein exchanged with Max Noether, who was then completing work on the lengthy report he and Alexander Brill published on the history of algebraic functions \cite{B-N}. This correspondence sheds new light on Noether's deep concerns about the importance of this report in substantiating his own priority rights and larger intellectual legacy.

\bigskip
\noindent
MSCCode: 14-03, 01A55, 01A60


\section{Introduction}

Julian Lowell Coolidge was a great expert on classical algebraic geometry \cite{Struik1955}. After studying under Corrado Segre in Turin, he went on to take his doctorate under Eduard Study in Bonn. As a mathematician, Coolidge excelled in writing books, some of them familiar to historians of geometry, others less so. Among the latter, his {\it Treatise on Algebraic Plane Curves} \cite{Coolidge1931} is easily overlooked. Certainly its style and contents appear very old-fashioned today, especially when set alongside a text like \cite{B-K}, despite obvious similarities in the subject matter. Coolidge cultivated an unusually informal writing style, even when describing  rather technical matters. He was also disarmingly honest, informing the reader whenever he happened to discuss a work without actually having held it in his hands.\footnote{Since he taught for many years at Harvard University, Coolidge had ready access to the rich holdings in Widener Library, whose director was his older brother \cite[670]{Struik1955}. So it did not often happen that he could not read a mathematical text firsthand.} 

One of the authors Coolidge  admired most deeply was Max Noether,\footnote{In his preface, Coolidge wrote: ``Large protions of the work are written according to the spirit and methods of the Italian geometers, to whom, indeed, the whole is dedicated [Ai Geometri Italiani, Morti, Viventi]. It would be quite impossible to describe the extent of the writer's obligation to them. Yet behind the Italians stands one whose contributions are even greater, Max Noether'' 
\cite[x]{Coolidge1931}.} whose publications from the 1870s exerted a lasting  influence on Italian algebraic geometers \cite{CES}. Noether followed in the footsteps of his principal mentor, Alfred Clebsch, who opened the doors to exploring the rich possibilities of Riemann's theory of complex functions for algebraic geometry, in particular the birational geometry of curves \cite[295--309]{Klein1926}.

The present paper discusses a central theorem in birational geometry first proved by Eugenio Bertini in \cite{Bertini1891}. Coolidge described the main ideas behind Bertini's proof in \cite[208--212]{Coolidge1931}, but he attributed the theorem to Clebsch. He did so owing to a short note that Felix Klein appended to the republication of Bertini's article in {\it Mathematische Annalen} \cite[160]{Bertini1894}. 
The precise circumstances that led to Klein's intervention have never been described before, but they 
can be easily reconstructed from letters Klein exchanged with Max Noether. This correspondence took place just as Noether and Alexander Brill were putting the last touches on their massive report on the history of algebraic 
functions \cite{B-N}.\footnote{After Noether's death, Brill clarified that they had adopted a clear division of labor in writing \cite{B-N}: Brill was responsible for the largely historical part up to and including Riemann's work, whereas Noether wrote about the various directions taken by contemporary researchers working in the wake of Riemann's novel innovations \cite{Brill1923}. In the discussion below, we will be concerned only with Noether's portion of the report. Noether chose to omit the more recent work on higher-dimensional algebraic varieties inaugurated by Italian researchers as well as the arithmetical approach taken by Dedekind and Weber. The latter as well as Hensel and Landsberg, {\it Theorie der algebraischen Funktionen einer Variabeln und ihre Anwendung auf algebraische Kurven und Abelsche Integrale} (1902) were taken up much later in a shorter report by Emmy Noether in \cite{Noether1919}. } 

This episode, as it emerges from the Klein--Noether correspondence, sheds new light on Noether's deep concerns about his report and the reception it might receive. Quite clearly, he saw parts of it as substantiating his own priority rights and larger intellectual legacy. Furthermore, the story told here testifies to the importance of oral communication not only for disseminating mathematical knowledge but also as a factor in contemporary discussions of priority claims. It can thus be seen as a case study supporting my
 longstanding interest in the oral dimensions of modern mathematical cultures (see  \cite{Rowe2003}, \cite{Rowe2004},  \cite{Rowe2018}).

\begin{figure}[htb]
\centering
\includegraphics[width=8cm,height=12cm]{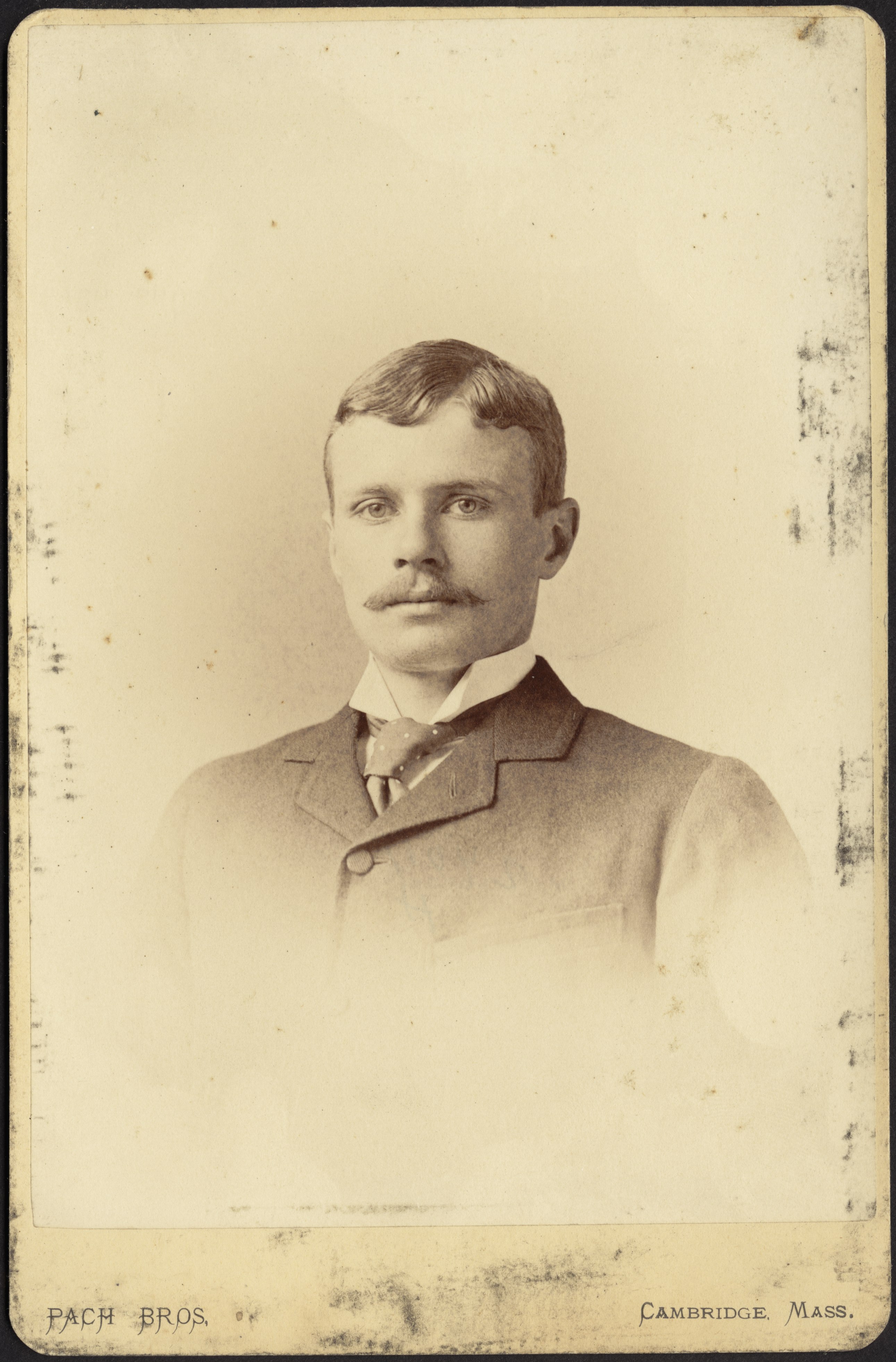} 
\caption{Julian Lowell Coolidge, 1873--1954
(Pach Bros., Cambridge, MA, photographer, ca. 1890)} 
\end{figure}

The result that Coolidge called ``Clebsch's Transformation Theorem'' states that any arbitrary plane algebraic curve can be birationally transformed into another having only double points as singularities.  Jean Dieudonn\'e, quite possibly on Coolidge's authority, also attributed this theorem to Clebsch  in \cite[37--38]{JD}, whereas Bertini described it as well known, part of the folklore of the era. He republished his paper \cite{Bertini1891} three years later in {\it Mathematische Annalen}, after noticing that prominent French authors -- among them Picard, Simart, and Poincar\'e -- continued to mention the theorem without citing his proof. On the other hand, Bertini failed to reference \cite{Noether1871} and other subsequent papers that showed how to transform a curve to obtain another with simple multiple points, i.e. distinct tangents for each branch. This oversight was perhaps entirely innocent, but in fact Noether's work represented the most difficult part in proving Clebsch's Theorem. 
As we shall see, this aspect helps to explain Noether's evident irritation in some of his letters to Klein.

\section{Klein's Correspondence with Noether} 

When Felix Klein received Bertini's note in early January 1894, what he read immediately set off alarm bells. So he turned to his long-time friend, Max Noether, asking him to help put out a potential fire.\footnote{Their friendship began in 1869 when both were studying under Clebsch in G\"ottingen; see \cite[40--48]{Tobies2019}.}
Klein quickly realized  that Bertini's proof of this folklore theorem in birational geometry was essentially the same as one he had learned about from Clebsch some 25 years earlier. He also vaguely  remembered writing to Noether about Clebsch's idea for removing higher-order singularities by exploiting the birational mapping of a plane to a cubic surface, as first presented in \cite{Clebsch1866}. Klein had only recently returned from his journey to the United States, during which he delivered his famous {\it Evanston Colloquium Lectures} \cite{Klein1894}. 

The first of these lectures dealt directly with the work of Clebsch and included these remarks about the above theorem:
\begin{quote}
Clebsch begins his whole investigation on the consideration of what he takes to be the most general type of an algebraic curve, and this {\it general curve} he assumes as having only double points, but no other singularities. To obtain a sure foundation for the theory, it must be proved that any algebraic curve can be transformed rationally into a curve having only double points. This proof was not given by Clebsch; it has since been supplied by his pupils and followers, but the demonstration is long and involved. \cite[4]{Klein1894}
\end{quote}
Here Klein cited the two papers \cite{B-N1874} and \cite{Noether1884}. Presumably, he had never studied these papers very carefully, so that when he read Bertini's note, which claimed to give the first real proof of this theorem, he realized that the assertion he had made in Evanston -- and that would soon be in print --  was mistaken.

Still, he made no mention of this in his letter to Noether from 4 January 1894.\footnote{Cod.Ms. Felix Klein, XII 637, SUB G\"ottingen.}
Instead, he gently suggested to Noether that he write Bertini, informing him that his proof was by no means new and that Clebsch had communicated it orally to Klein long ago. In the letters that follow it should be borne in mind that Klein had been the principal editor of {\it Mathematische Annalen} for nearly twenty years, whereas Noether only joined the board as an associate editor in 1893.
\begin{quote}
You will perhaps want to attach a comment from your side to the note by Bertini included here. I send it to you with the request that, if necessary, you take up correspondence with the author and also pass the note on to Dyck [Walther Dyck was the managing editor of {\it Mathematische Annalen}]. The offprint itself I would like to have returned because I took it from my own copy of the Revista.\footnote{Bertini presumably had sent Klein an offprint of \cite{Bertini1891} along with a request that this short note be reprinted in the {\it Annalen}.} As far is the matter itself is concerned, I remember that Kronecker and Clebsch discussed this during the fall of 1869 in Berlin, which led me to the basic method I've since often explained and which is essentially the same as Bertini's: I view the plane as the image of a cubic surface and then project the curve carried onto the surface back again into a plane, etc. etc.

In any case, I wrote or spoke to you about this once. If you, as I suspect, wish to add a note to Bertini's article, perhaps you could take the opportunity to say a word about this as well.\footnote{Beifolgende Notiz von Bertini wirst Du vielleicht mit einer Bemerkung von Deiner Seite versehen wollen. Ich schicke sie daher zu, mit der Bitte, wenn n\"otig mit dem Verf[asser]
selbst zu correspondieren und \"ubrigens die Note an Dyck. weitergehen zu lassen. Den Druckbogen selbst m\"ochte ich mir sp\"ater zur\"uck erbitten, da ich ihn aus meinem Exemplar der Rivista ausgel\"ost habe.  Was die Sache angeht, so erinnere ich mich, dass Kronecker und Clebsch im Herbst 1869 in Berlin dar\"uber verhandelten, worauf ich den Ansatz fand, den ich seither oft vortrug und der im Wesen mit dem von Bertini identisch ist: ich sah die Ebene als ``Abbildung'' einer Fl\"ache dritter Ordnung an und projecierte dann wieder die auf diese Fl\"ache \"ubertragenen Curve auf eine Ebene, etc. etc. 

Ich habe Dir jedenfalls einmal davon geschrieben oder gespro\-chen. Wenn Du, wie ich vermuthe, von Dir aus [eine] Note zu dem Aufsatz von Bertini hinzuf\"ugen willst, nimmst Du vielleicht Gelegenheit auch hier\"uber ein Wort einfliessen zu lassen.}
\end{quote}

Noether was out of town when Klein's letter arrived, but on returning to Erlangen he wrote back in an agitated state of mind. At first, he could not believe that  Klein had actually written to him about this matter long ago, and that neither of them had ever spoken about it since. 
After looking through Klein's letters from late 1869, however, he found that this was, indeed, the case. At that time, Klein had been studying in Berlin along with his new-found Norwegian friend, Sophus Lie, while Noether was working under Clebsch in G\"ottingen.  During the fall vacation, Clebsch visited Berlin for a few days, during which time he and Kronecker discussed methods for desingularizing an algebraic curve. Klein then learned about this discussion from Clebsch and wrote the following in a letter to Noether from 17 December 1869:\footnote{Cod.Ms. Felix Klein, XII 527, SUB G\"ottingen.}
\begin{quote}

You will perhaps also be interested in what I learned from Clebsch when he was here during the fall vacation. (I really don't know any more whether I already wrote you this or not.) Kronecker has  proved, namely, that a plane curve with arbitrary singularities can always be transformed into one with a single multiple point whose tangents are all distinct. Clebsch then pointed out that one can easily resolve this multiple point into ordinary double points. One views the plane containing the curve as the image of an $F_3$ [cubic surface] in such a way that a fundamental point passes into the given multiple point. -- This comment of Clebsch appears to me to possess considerable mathematical value. If I'm not mistaken, it allows for a type of extension so that one can immediately reduce arbitrary singularities to ordinary double points.\footnote{Dann wird Dich vielleicht noch interessieren, was ich von Clebsch erfahren habe, als er in den Herbstferien hier war. (Ich wei{\ss} wirklich nicht mehr, ob ich es Dir nicht schon geschrieben habe.) Kronecker n\"amlich hat nachgewiesen, da{\ss} sich eine ebene Kurve mit beliebigen Singularit\"aten immer auf eine solche zur\"uckf\"uhren  l\"a{\ss}t, die einen einzigen mehrfachen Punkt besitzt, dessen Tangenten s\"amtlich verschieden sind. Damals machte Clebsch darauf aufmerksam, da{\ss} man diesen Punkt nun sehr einfach in gew\"ohnliche Doppelpunkte aufl\"osen kann, indem man die Ebene der Kurve als Bild einer $F_3$ auffa{\ss}t, wobei ein Fundamentalpunkt 
in den gegebenen mehrfachen Punkt r\"uckt. -- Mir scheint diese Bemerkung von Clebsch einen hohen mathematischen Wert zu besitzen. Irre ich nicht, so l\"a{\ss}t sie sich in einer solchen Art erweitern, da{\ss} man beliebige Singularit\"aten unmittelbar auf gew\"ohnliche Doppelpunkte reduzieren kann.}

\end{quote}

One can easily imagine Noether's astonishment when he read this passage, which not only clearly indicated that Clebsch had found a simple method for reducing the singularities of a plane curve to double points. It also indicated that Kronecker may have  anticipated Noether's own argument using quadratic transformations to show how to obtain a curve having only ordinary singularities -- those whose branches have distinct tangents.

 Since Noether was very familiar with the mapping in \cite{Clebsch1866}, he would have realized immediately how it can be used to desingularize a given mutliple point. Clebsch's mapping takes  six fundamental points in general position in the plane and blows these up into six of the 27 lines on a cubic surface $F_3$ (Fig. 2). Each set of five fundamental points determines a conic, and these six conics also go over into six lines. The remaining 15 are the images of the lines connecting pairs of the six fundamental points. If one then starts with a plane algebraic curve $C_n$ with singular point $P\in C_n$ of multiplicity $m$, then by letting $P$ coincide with one of the six fundamental points of a Clebsch mapping has the effect of blowing up this singularity into $m$ points on the image line. The curve $C_n$ will then pass over to a space curve $C'$  lying on $F_3$. The remainder of the argument then  involves carefully finding a point $Q\in F_3$ so that the projection of $C'$  back into the plane will only introduce double points as singularities. As Coolidge shows in some detail \cite[210--212]{Coolidge1931}, the point $Q$ has to be chosen so as not to fall on a line that happens to belong to any of four different 1-parameter systems, which of course can always be done. The image of $C'$  is then a new planar curve $C''$ on which the singular point $P\in C_n$  now corresponds to $m$ nonsingular points; moreover, the only new singularities that arise will be simple nodes. Repeating this argument for each higher singularity then yields a new curve birationally equivalent to $C_n$ and whose only singularities are double points.

\begin{figure}[htb]
\centering
\includegraphics[width=12cm,height=7cm]{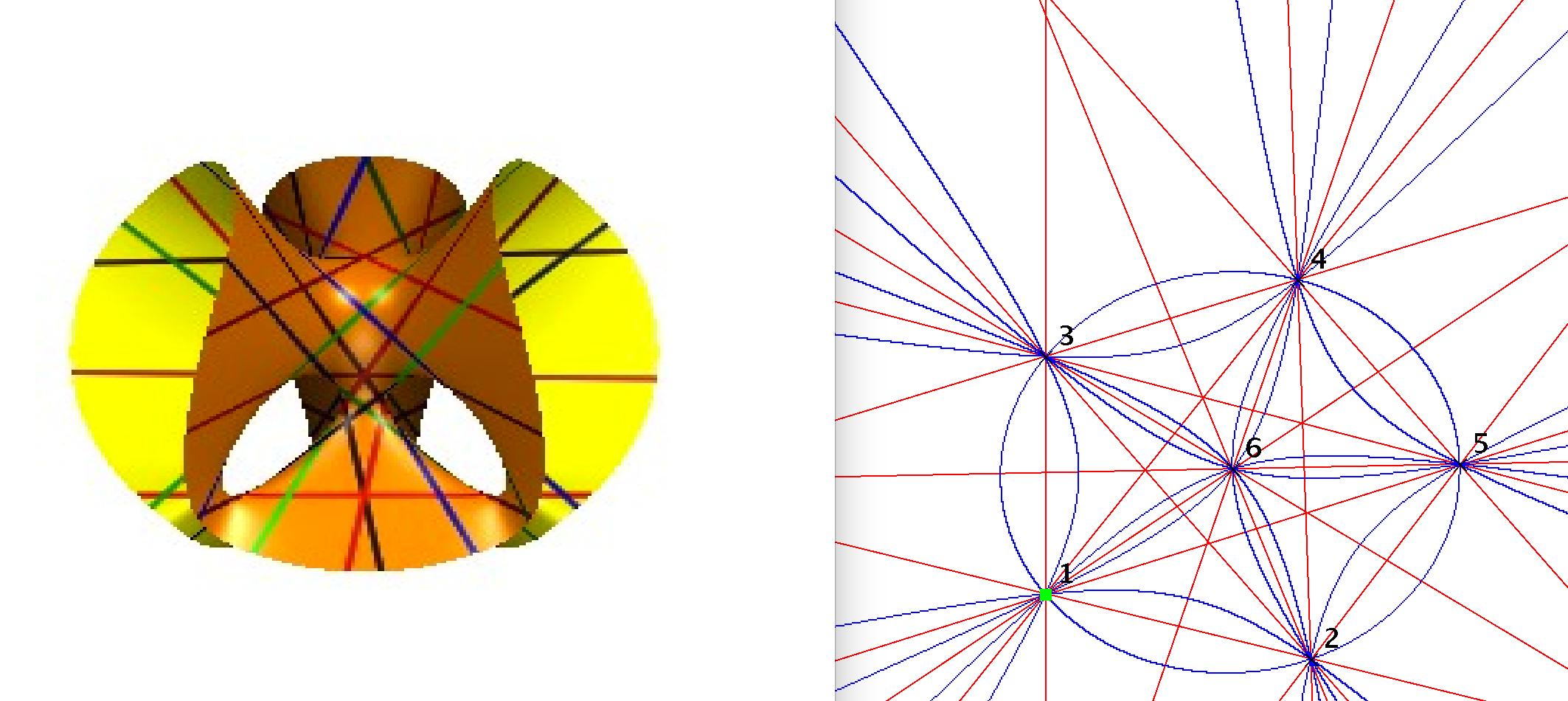} 
\caption{Clebsch's Mapping of a Plane to a Cubic Surface from \cite{Clebsch1866}
(Pictures Courtesy of Oliver Labs)} 
\end{figure}

Klein's request that Noether inform Bertini of these oral communications from 1869 left Noether at a loss as to what he should  do. He regarded  Klein's intervention as a very belated  and also 
highly unexpected  attempt to set the record straight.  Under normal circumstances, Noether would have surely wanted to clarify such matters as well. This particular case, however,  put  him in a very awkward situation. He was certainly less concerned about  informing Bertini regarding these matters. What worried him about the revelations in Klein's letter from 1869 was that this new information might have deeper implications for assessing priority claims with regard to the resolution of singularities. If so, this would presumably require that Noether reconsider
 what he 
 had recently  written about this topic for the forthcoming report \cite{B-N}. That thought left him deeply troubled. Since Noether's own work on this topic had played a central role, he clearly had a  personal   interest in clarifying the situation. In particular, he wanted to ensure that what he wrote in no way  clashed with Klein's views on these matters. His concern about a potential conflict of this sort thus prompted him to write a lengthy and somewhat defensive letter to Klein on 12 January.\footnote{Cod.Ms. Felix Klein, XI 112, SUB G\"ottingen.}
\begin{quote}
Due to my absence from Erlangen nothing was done with the mailing from the $4^{th}$ regarding Bertini. I also cannot take care of it without first reaching an understanding with you.

Your message, that you had written to me in \underline{1869} that you had resolved a singular point by repeated use of plane mappings of cubic surfaces, astonished me incredibly. I have absolutely no recollection of this, and also know precisely, that my note for the G\"ottingen Nachrichten
\cite{Noether1871}  in \underline{June 1871} was completely independent, and I cannot understand why in 1871 neither you nor Clebsch ever mentioned that alleged communication from 1869 or why you have since then, until now, never come back to this, for example in the book Clebsch-Lindemann.\footnote{Durch meine Abwesenheit von Erlangen ist diese, Bertini betreffende Sendung vom $4^{ten}$ liegen geblieben. Ich kann dieselbe auch nicht erledigen, ohne mich vorher mit Dir verst\"andigt zu haben. 

Deine Mitteilung, da{\ss} Du im Jahre \underline{1869} die Aufl\"osung eines singul\"aren Punktes durch wiederholte Benutzung der eb[enen] Abbildungen von Fl\"achen $3^{ter}$ O[rdnung] bewerkstelligt hattest und mir dies damals auch geschrieben hattest, hat mich auf's \"Au{\ss}ersten frappiert. Ich erinnerte mich absolut nicht daran, wei{\ss} auch v\"ollig genau, da{\ss} ich im \underline{Juni 1871} beim Verfassen meiner Note  in den G\"ottinger Nachrichten  [dieser?] durchaus unabh\"angig war, und kann nicht begreifen, warum [?] 1871 weder Clebsch noch Du mit irgend einem Wort mir gegen\"uber jemals auf jene angebliche Mitteilung von 1869 hinwiesen oder warum Du auch seitdem, bis jetzt, nicht darauf
 zur\"uckgekommen bist, z.B. nicht im Buche Clebsch-Lindemann.}
\end{quote}

Noether's irritation over having been left in the dark was evident, but he saw from the passage in Klein's old letter that Clebsch and Kronecker had, indeed, discussed this matter. Furthermore, Klein had already contemplated using the Clebsch 
mapping as a method for reducing singularities to double points. He, therefore, felt compelled to inform Klein about what he had written regarding Kronecker's contributions as well as his own for the report \cite{B-N}.

\begin{quote}
I now went through your correspondence from 1869, which confirms your message, as shown by the enclosed letter. So now I have to deal with this not only due to Bertini's note but also owing to my report on algebraic functions.\footnote{Deine Correspondenz von 1869, die ich nun durchsuchte, best\"atigt Deine Mitteilung, wie Dein hier mitfolgende Brief zeigt. Ich habe mich nun mit Dir nicht nur wegen der Bertini'schen Note, sondern auch wegen meines Referates \"uber algebr. Funktionen, auseinanderzusetzen.}
\end{quote}
Noether then proceeded to summarize the relevant parts of his report, in some places citing verbatim from the text itself. He then continued with a provisional assessment of the various discoveries from 1869:

\begin{quote}
I wrote about all this in such detail because I must know whether you wish to make a priority claim against any part of this [report].

So far as I can see, 1) Kronecker was incorrect, if he opined, that the curve in his article on the discriminant has only \underline{one} \underline{ordinary} multiple point; 2) that Clebsch was the first to use the mapping of an $F_3$ for the resolution of  \underline{one} ordinary point; that you used the $F_3$ successively  to resolve several ordinary singular points; 4) you do not show that this method is effective for arbitrary singular points; 5) nor, in particular, that the process terminates; 6) you give no applications of it.\footnote{Ich habe dies Alles so ausf\"uhrlich geschrieben, weil ich wissen mu{\ss}, ob Du gegen irgend einen Teil desselben 
eine Priorit\"atsreclamation geltend zu machen w\"unschst. 

So viel ich sehe, hat 1) Kronecker nicht Recht, wenn er meint, da{\ss} in seinem Discr[iminaten] Aufsatz seine  [?] Curve nur \underline{einen} \underline{gew\"ohnlichen}  mehrfachen Punkt hat; 2) hat Clebsch zuerst die Abbildung der $F_3$ zur Aufl\"osung 
\underline{eines} gew\"ohnlichen Punktes benutzt; 3) hast Du die succ[essive] $F_3$ zur Aufl\"osung  mehrerer gew\"ohnlichen sing. Punkte benutzt; 4) f\"uhrst Du nicht aus, da{\ss} diese Methode bei beliebigen  singul\"aren Punkten in der That [wirksam?] ist; 5) zeigst Du insbesondere nicht, da{\ss} der Proze{\ss} abschlie{\ss}t; 6) gibst Du keine Anwendung davon.}
\end{quote}

At this point, Noether briefly listed some of  his own accomplishments, obviously in order to make clear why he saw no need to pay great heed to the problem of reducing ordinary singularities to double points. ``I did not give any further resolution of multiple points because I held this to be unnecessary, and above all, because this seemed to then as it does now completely evident.'' (Die weitere Aufl\"osung der vielf[achen] Pktn habe ich nicht mitgegeben, weil ich sie f\"ur unn\"otig hielt, und vor Allem: weil sie mir v\"ollig evident schien und noch scheint.'') Nevertheless, he asked Klein to state  his views about these matters as clearly as possible:

\begin{quote}
If your priority claim is directed toward the last part V [of the report] (ordinary multiple points) -- which alone concerns Bertini etc., I have nothing against adding some  information if that is suitable for you; if you however today believe, after such a long time and after your silence throughout the year 1870 when we corresponded with one another, 
 that you have grounds for a priority claim for the resolution of 
\underline{singular} points based on a remark that passed over me in 1869 without a trace, then that would not pertain to Bertini but rather to me: I would appear in rather false light as having known of Kronecker's [result] and your successive method without citing these. I cannot see any good way to reach a compromise in this last respect, and so I must first ask you in this second, for me completely unanticipated case for precise indications as to how far you are now in disagreement with the above contents from my report.\footnote{Richtet sich nun Dein Priorit\"atsanspruch auf den letzteren Teil V (gew\"ohnliche vielfache Punkte) -- der bei Bertini etc. allein in Betracht kommt, so habe ich nichts dagegen, eine Zuf\"ugung zu machen, wenn es Dir pa{\ss}t; wenn Du aber glaubst, auch f\"ur die Aufl\"osung der \underline{singul\"aren} Punkte heut, nach so langer Zeit, und nachdem Du im Jahr 1870, wo wir doch in Correspondenz standen, geschwiegen, einen Anspruch begr\"unden zu k\"onnen -- auf ein Aper\c{c}u hin -- das an mir 1869 spurlos vor\"uberschwand -- so w\"urde sich das nicht auf Bertini, sondern auf mich beziehen: ich k\"ame in eher falsches Licht, als h\"atte ich Kronecker und Dein succ[essive] Verfahren gekannt, ohne es zu nennen. Ich sehe in dieser letzteren Beziehung keinen richtigen Weg des Ausgleiches, und m\"u{\ss}t  Dich in diesem zweiten, mir g\"anzlich unvermuthet kommenden Falle erst um genaue Angabe dar\"uber bitten, wie weit Du nun mit dem obigen Inhalte meines Referates nicht einverstanden sein wirst.}
\end{quote}

As for Klein's original request -- that Noether consider writing a note that would be appended to Bertini's reprinted article, mentioning the state of affairs in 1869 -- he saw no reason why he should add any comment whatsoever:

\begin{quote}
Furthermore, I remark that I have no real motivation to comment on Bertini's note, since the question concerning the transformation of a curve with ordinary multiple points into another with double points is a matter of complete indifference to me -- and evident, as I already said. That Bertini neglected to cite me in connection with the passage from a singular curve to one with multiple points\footnote{Bertini began his proof with a single sentence asserting that one should first transform the curve by a suitable Cremona transformation to obtain another having only ordinary multiple points  \cite[159]{Bertini1894}.} should hardly matter; that is of course known. Probably I will send the note on to Dyck without any comment.\footnote{Noch bemerke ich, da{\ss} ich zur Bertini'schen Note eine Anmerkung zu machen eigentlich gar kein Anla{\ss} habe, da mir die Frage der \"Uberf\"uhrung einer Curve mit gew[\"ohnlichen] mehrfachen Punkten in eine solche mit Doppelpunkten ganz gleichg\"ultig ist -- und evident, wie ich schon sagte. Da{\ss} B[ertini] mich in der Frage der \"Uberf\"uhrung einer sing[ul\"aren] Curve auf eine solche mit mehrfachen Punkten nicht citiert, soll wohl nichts hei{\ss}en; das ist ja bekannt. Wahrscheinlich w\"urde ich die Note unvermerkt an Dyck weiter schicken.}
\end{quote}

Klein was clearly more than a little surprised when he read 
Noether's long response to his original inquiry. This response also included  Klein's own letter from 17 December 1869 containing  the key passage cited above. He probably never imagined that Noether would consider writing about these oral communications in \cite{B-N}, but in fact this official report for the German Mathematical Society went well beyond the literature found in journals and books. In particular, Noether wrote at length about Kronecker's unpublished program for desingularizing algebraic curves, an approach he contrasted with his own.\footnote{In \cite[370]{B-N} Noether cited Kronecker's statement that he had communicated his methodological views to Riemann and Weierstrass in 1858, presented these ideas to the Berlin Academy in 1862, and expounded the same in his lectures beginning in the winter semester of 1870/71. They were only published in  \cite{Kronecker1881}, which announces a sequel that never appeared.}  His report also addressed in great detail the contributions of Karl Weierstrass, which required many references to unpublished results from  the latter's  lecture courses. Since the 1870s, Noether and Klein had both studied Weierstrass' work very avidly by means of various {\it Ausarbeitungen} made by his students in Berlin. In many instances, these served as the basis for various priority claims made by Weierstrass, but especially by his closer associates.

 Noether thus had every reason to take Klein's overlooked letter from 1869 quite seriously. Surely he felt very relieved to learn that Klein  had no intention at all of raising a priority claim. Probably Klein never imagined that Noether should consider revising his report in order to bring out these early discussions of  ``Clebsch's Theorem,'' but since Noether raised this possibility himself, Klein now made a simple suggestion in this direction:\footnote{Klein to Noether, 19 January 1894, Cod.Ms. Felix Klein, XII 638, SUB G\"ottingen.}

\begin{quote}
This is truly a very strange situation. I did not want to direct any priority claim against you recently, all the less so as I had myself completely forgotten the contents of my letter, insofar as it stood in competition with your works. My intention was only to take a stand against Bertini, and I also only do that because I would otherwise appear in false light. The point is namely this, that often in recent years -- and, in particular, also recently in the New York Mathematical Society, where I was asked about the resolution of multiple points with distinct branches -- I have lectured about this, naturally without any mention of Bertini. If now immediately afterward I publish Bertini's note in the Annalen it creates the impression that I had intentionally neglected to mention Bertini's name. I thus believe it is necessary that I add a comment to Bertini's note, which I include for provisional passage to Dyck. That way I will create no problems for you.

On the other hand it would appear to me correct if you would perhaps include in your report the passage from my letter \dots with the comment that we both only later became aware of this communication, and that I naturally do not want to make any priority claim based on a sketch as opposed to a thorough investigation of the object in question. But I leave that entirely up to you, and if you decide not to do ths I will certainly not raise any objection.\footnote{Das ist wirklich eine sehr merkw\"urdige Situation. Ich habe neulich gegen Dich gar keine Priorit\"atsreclamation richten wollen, um so weniger, als ich den Inhalt meines Briefes, soweit er mit Deinen Arbeiten in Concurrenz tritt, selber vollkommen vergessen hatte. Meine Absicht war nur gegen Bertini Stellung zu nehmen. Und auch dies thue ich nur, weil ich sonst in schiefes Licht komme. Die Sache ist n\"amlich die, dass ich den letzten Jahren \"ofter und insbesondere auch neulich in der New Yorker Mathematical Society, als ich nach der Aufl\"osung der vielfachen Puncte mit getrennten Aesten gefragt wurde, die Sache vorgetragen habe, nat\"urlich ohne Bertini zu nennen. Wenn ich nun unmittelbar hernach die Note von Bertini in den Annalen drucke, so entsteht der Eindruck, als habe ich Bertini's Namen absichtlich verschwiegen. Ich glaube also, dass es nothwendig ist, dass ich der Note von Bertini einen Zusatz mache, wie ich ihn zu ev. Weiterbef\"orderung an Dyck beilege. Damit gerathe ich Dir gar nicht in's Gehege. 

Andererseits schien mir richtig, dass Du vielleicht in Deinem Referate die Stelle meines Briefes zwischen den beiden jetzt von mir am Rande angebrachten Sternchen unter der Seite abdrucktest, mit dem Bemerken, dass wir beide erst hinterher wieder auf diese Mitteilung aufmerksam geworden sind und dass ich selbst\-verst\"andlich auf ein Aper\c{c}u keinen Priorit\"atsanspruch gegen\"uber einer ausf\"uhrlichen Durcharbeitung des Gegenstandes gr\"unden will. Aber ich \"uberlasse Dir das vollkommen  und werde gewi{\ss} nicht, wenn Du es nicht thust, reclamieren.}
\end{quote}

Klein surely found Noether's letter a quite bizarre and certainly hypersensitive overreaction to this whole matter. Probably he remembered nothing about the details of the discussion between Kronecker and Clebsch. On the other hand, Noether's lengthy remarks to Klein 
 about  Kronecker's methods and claims -- remarks based on the not yet published text of \cite{B-N} -- reveal how carefully he  approached this terrain. Did Klein even bother to read these parts of Noether's letter? He was a busy man; futhermore,  he fully accepted Noether's authority as {\it the} leading expert on all such matters.

\section{Klein's Note and Noether's Commentary}

Noether's relationship with Klein had always been harmonious, and it would remain that way in the future, despite occasional  differences with respect to {\it Mathematische Annalen}, the journal co-founded by their mutual mentor, Alfred Clebsch. Once he saw that Klein merely wanted to clarify his own early and certainly quite  marginal involvement with what Noether regarded as a theorem of no great importance for the resolution of singularities  for curves, the latter  was more than  happy to accomodate him.
Klein's appended  {\it Zusatz}  to \cite{Bertini1894} reads:

\begin{quote}
The method of Bertini, speaking geometrically, amounts to regarding the plane, in which the curve with a singular point lies, as a single-valued image of a cubic surface, by which the curve is transformed into a space curve, and the latter is projected into another plane by means of a sufficiently general [projection] point. In this formulation I am familiar with the Ansatz through Clebsch, who communicated it to me orally in the fall of 1869. I mention this only because I have especially in recent times often appealed to this in my courses and lectures. The readers of the Annalen will have no less reason to be thankful to Mr. Bertini for his detailed presentation.\footnote{Die Methode von Bertini kommt geometrisch zu reden darauf zur\"uck, die Ebene, in welcher uns die Curve mit singul\"arem Punkte  gegeben ist, als eindeutige Abbildung einer Fl\"ache dritter Ordnung zu betrachten, dadurch die Curve in eine Raumcurve zu verwandeln und letztere hinterher wieder von einem hinreichend allgemeinen Punkt aus auf eine andere Ebene zu projiciren. In dieser Form ist mir der Ansatz noch von Clebsch her bekannt, der mir denselben im Herbst  1869  m\"undlich mittheilte. Ich erw\"ahne dies nur, weil ich gerade in letzter Zeit in meinen  Vorlesungen und Vortr\"agen wiederholt daran angekn\"upft hatte. Die Leser der Annalen werden darum Hrn. Bertini f\"ur seine ausf\"uhrliche Darstellung nicht geringeren  Dank wissen.}
\end{quote}

This was the note that Coolidge mentioned in attributing  to Clebsch the theorem stating that every algebraic curve can be birationally transformed to another having only double points as singularities  \cite[212]{Coolidge1931}. Oddly enough, Coolidge seems to have overlooked the nearly contemporaneous note that Noether, acting on Klein's suggestion, added to his report. This contains the main substance of the passage from Klein's letter from 1869, as given above. Noether placed this as a footnote to a sentence stating that \cite{Noether1871} gave the first (incomplete) proof  for resolving the singularities of an algebraic curve \cite[371]{B-N}. The footnote begins by downplaying
 the critical passage, ``which at the time remained {\it completely} ignored by [Noether] and which was just recently brought to his and the author's attention by accident''  (``welche damals von [Noether] {\it g\"anzlich} unbeachtet geblieben war, und auf welche seine, wie  des Schreibers Aufmerksamkeit erst jetzt wieder zuf\"allig  gelenkt wurde'').  Noether might have left the matter there, but he decided to add some commentary to be sure that nothing in Klein's letter left the reader wondering.

\begin{quote}
Clebsch's observation only concerns a special form of his general investigations on the transformation of an $h$-fold point with distinct tangents into $h$ simply separated points \dots ; the significance of Klein's observation, which was given only as a mere sketch without further elaboration, cannot be judged. In any case it was not exploited for the resolution of truly {\it singular} points. It could thus not serve as the basis for a priority claim.\footnote{Die Bemerkung von Clebsch bezieht sich nur auf eine specielle Form seiner allgemeinen Betrachungen \"uber Transformation eines $h$-fachen Punktes mit getrennten Tangenten in $h$ einfach getrennte Punkte \dots ; die Tragweite der Klein'schen Bemerkung l\"asst sich, wie dieselbe als blosses Aper\c{c}u ohne weitere Ausf\"uhrung gegeben ist, nicht beurteilen, jedenfalls ist sie f\"ur die Aufl\"osung der wirklich {\it singul\"aren} Punkte nicht verwertet worden. Einen Priorit\"atsanspruch k\"onnte dieselbe also nicht begr\"unden.}
\end{quote}

Whether or not Coolidge ever read these remarks, he certainly did not share Noether's dismissive attitude with regard to ``Clebsch's Theorem,'' based on the argument in \cite{Bertini1894}. This is apparent not only from his lengthy exposition of it, but also from some remarks he made in the concluding section of \cite{Coolidge1931}. There he describes various systems of plane curves that remain invariant under finite groups of Cremona transformations. In this connection, he called attention to a paper by Anders Wiman and wrote: ``In many cases there is much to be gained by reverting to the method of mapping the plane on an auxiliary cubic surface that we developed \dots to prove Clebsch's transformation theorem \dots . The problem then becomes one of finding groups of collineations of a three-dimensional space which  leave a certain cubic surface invariant'' \cite[498]{Coolidge1931}. The paper Coolidge cited,  \cite{Wiman1896}, appeared just two years after Noether's report and was published, appropriately enough, in Clebsch's journal, {\it Mathematische Annalen}.

\end{document}